\documentstyle[12pt]{article}
\newlength{\abstractwidth}
\renewcommand{\thefootnote}{\fnsymbol{footnote}}
\renewcommand{\thanks}[1]{\footnote{#1}} % Use this for footnotes
\newcommand{\starttext}{
\setcounter{footnote}{0}
\renewcommand{\thefootnote}{\arabic{footnote}}}

\newcommand{\be}{\begin{equation}}
\newcommand{\bea}{\begin{eqnarray}}
\newcommand{\eea}{\end{eqnarray}}
\newcommand{\ee}{\end{equation}}

\def\ba{\begin{eqnarray}}
\def\ea{\end{eqnarray}}

\def\v{\vskip .1in}

\def\al{\alpha}

\def\m{\mu}

\def\o{\omega}

\def\r{\rho}

\def\O{\Omega}
\def\T{\Theta}

\def\cB{{\cal B}}

\def\cH{{\cal H}}

\def\cL{{\cal L}}
\def\cM{{\cal M}}

\def\cO{{\cal O}}

\def\cX{{\cal X}}

\def\R{{\bf R}}
\def\C{{\bf C}}
\def\P{{\bf P}}

\def\Aut{{\rm Aut}}

\def\ti\tilde

\def\u{\underline}

\def\pl{\partial}

\def\i{\infty}

\def\s{\sum}

\def\ddb{\partial\bar\partial}
\def\sub{\subseteq}
\def\ra{\rightarrow}

\def\ti{\tilde}

\def\o{\omega}

\def\[{{\bf [}}
\def\]{{\bf ]}}

\def\pl{\partial}

\v

\begin{document} \starttext \baselineskip=15pt
\setcounter{footnote}{0} \newtheorem{theorem}{Theorem}
\newtheorem{proposition}{Proposition}
\newtheorem{lemma}{Lemma} \newtheorem{corollary}{Corollary}
\newtheorem{definition}{Definition} \begin{center} {\Large \bf ON THE
REGULARITY OF GEODESIC RAYS ASSOCIATED TO TEST CONFIGURATIONS}
\footnote{Research
supported in part by National Science Foundation grants DMS-02-45371
and DMS-05-14003}
\\
\bigskip

{\large D.H. Phong$^*$ and
Jacob Sturm$^{\dagger}$} \\

\bigskip

\end{center}
\v
\v\v\v

\begin{abstract}

Geodesic rays of class $C^{1,1}$ are constructed for any test configuration
of a positive line bundle $L\to X$, using resolution of singularities.
The construction reduces to finding a subsolution of the corresponding
Monge-Amp\`ere equation. Geometrically, this is accomplished by the use
a positive line bundle on the resolution which is trivial outside
of the exceptional divisor.

\end{abstract}

\bigskip
\baselineskip=15pt
\setcounter{equation}{0}
\setcounter{footnote}{0}

\section{Introduction}
\setcounter{equation}{0}

Let $X$ be a compact complex manifold, $L\ra X$ an ample line
bundle, and $\cH$ the space of K\"ahler metrics $\o$ on $X$ with
$\o\in c_1(L)$. The work of Mabuchi \cite{M87}, Semmes \cite{S92}
and Donaldson \cite{D99} shows that $\cH$ is a infinite-dimensional
non-positively curved symmetric space with respect to a natural
Riemannian structure. Moreover, according to the
program proposed by Donaldson \cite{D99},
the existence of metrics
of constant scalar curvature in $\cH$ is linked to the existence
of geodesic rays in $\cH$, and the uniqueness of such metrics
is linked to the existence of geodesic segments in $\cH$ (which
we shall refer to as Donaldson rays and Donaldson segments).
On the other hand, the conjecture of Yau-Tian-Donaldson \cite{Y, T97, D02}
gives a necessary and sufficient criterion for the existence
of constant scalar curvature metrics which can be formulated in terms of Bergman
geodesics.
(Recall that the space $\cH_k$ of Bergman metrics is just the space
of pullbacks of the Fubini-Study metric by the Kodaira imbeddings defined by bases of
$H^0(X,L^k)$.)
Thus it is natural to look for a direct link
between Bergman segments/rays and Donaldson segments/rays.

\v

This program was carried out for geodesics segments in \cite{PS06}: If
$h_0,h_1\in\cH$, then we showed that  there is a sequence of
Bergman geodesics $h_t(k)$, with $0\leq t\leq 1$,  whose limit is a $C^{1,1}$
Donaldson geodesic $h_t$ joining $h_0$ to $h_1$. The
proof of this theorem consists of three main steps:
\v
1) The limit of the $h_t(k)$ is a weak geodesic $h_t$ joining $h_0$ to $h_1$.

2) There exists some $C^{1,1}$ geodesic $\ti h_t$ joining $h_0$ to $h_1$.

3) The two geodesics are equal: $h_t=\ti h_t$.

\v
The first step was carried out in \cite{PS06}.
The main ingredients  are the Tian-Yau-Zelditch \cite{Z} expansion and
the monotonicity theorem of Bedford-Taylor \cite{BT76}.
\v
The second step is due to  Chen \cite{Ch}. The main ingredients are the interior
$C^2$ estimates of Yau \cite{Y78}, the general theory of Caffarelli, Kohn, Nirenberg,
and Spruck \cite{CKNS}, and the
boundary $C^2$ estimate of Guan \cite{G} for the  Monge-Amp\`ere equation.
The last part makes use of a blow-up
analysis of the solution.
\v
The third step is also in \cite{PS06}. The main ingredient is the pluri-potential capacity
theory of Bedford-Taylor \cite{BT82}.

\v

Now we turn to the program for geodesic rays. We would like
to implement the same three steps which were used to solve
the geodesic segment problem. To describe the results, we first
let $\ti\cH$ be the space of
positively curved metrics on $L$, so that $\ti\cH/\R$
can be identified with $\cH$
(it is somewhat easier to formulate the results for
$\ti\cH$). Step 1) for geodesics rays was carried out
in \cite{PS06a}. We showed that if $h_0\in \ti\cH$, and
if $\r:\C^\times\ra \Aut(\cL\ra\cX\ra\C)$ is
a test configuration $T$, then there  exists
a weak geodesic ray emanating from $h_0$ and
``pointing in the direction'' of~$T$. This
ray is defined as a limit of a certain natural sequence of
Bergman rays $h_t(k)$
with $t\in [0,\i)$.  The
proof uses again the Tian-Yau-Zelditch expansion
and the Bedford-Taylor monotonicity theorem. The
main new ingredients are Donaldson's imbedding
theorem for test configurations (which is spelled
out in detail in \cite{PS06a}) and  a careful analysis of
the asymptotic expansion of the Bergman rays as $t\ra\i$.

\v

The purpose of this note is to carry out  step 2) for
geodesic rays.
We shall show that if $h_0\in\ti\cH$ and
if $T$ is a test configuration, then there is a
$C^{1,1}$ geodesic ray $\ti h_t$ emanating from $h_0$
with the following property: $h_t$ extends to a solution
of the Dirichlet problem
for the Monge-Amp\`ere equation on a resolution of
singularities $\ti\cX\ra\cX$.
To do this, we first construct a subsolution of the
associated Monge-Amp\`ere equation. The existence
of a subsolution is not automatic, since the underlying
manifolds are not strictly pseudo-convex. After
the subsolution is constructed, the arguments of \cite{CKNS, G, Ch}
are used
to construct the geodesic
$\ti h_t$.

\v

In order to complete the program for geodesics rays,
we need to implement step 3),
which is the identification of
$h_t$ with $\ti h_t$. This will be considered in a forthcoming
paper.
\v

Recently, there has been considerable activity in the study of
Donaldson geodesics. In \cite{AT}, Arezzo-Tian have constructed
families of analytic geodesic rays near infinity for test configurations
with smooth central fiber. In \cite{Ch06}, Chen has constructed
geodesic rays parallel to a given one, under various assumptions.
In \cite{SZ}
Song and Zelditch have provided a sharp analysis of the
approximation in \cite{PS06} of Donaldson geodesics by Bergman
geodesics, for toric varieties. In particular, prior to our present work, they \cite{SZf} have
shown that the optimal regularity for {\it geodesic rays}
constructed in this manner for toric varieties in general is
$C^{1,1}$.

\section{Statement of theorem}

Let $X$ be a compact complex
manifold and $L\ra X$ an ample
line bundle.
Recall that a test configuration $T$ for $(X,L)$
is a homomorphism
\be
\label{rho}
\r:\C^\times\ra\Aut(\cL\ra\cX\ra\C)
\ee
where
$\pi:\cX\ra\C$ is a proper flat morphism, $\cL\ra\cX$
is a line bundle which is ample on the fibers and
$\r$ is an action of $\C^\times$ on $\cL\ra\cX\ra\C$
which covers the standard action of $\C^\times$ on $\C$,
and which satisfies the following additional property: $\cL_1=L$ and $X_1=X$
(where, for $w\in\C$,
we let $X_w=\pi^{-1}(w)$ and $L_w=\cL|_{X_w}$). We shall
say $T$ is trivial if $\cX=X\times \C$ and $\r$ is the
trival homomorphism.
\v
Fix a resolution of singularities
\be
p:\ti\cX\ra\cX\ra \C.
\ee
According to \cite{BM, V} (see also \cite{H}), the resolution
can be chosen to be equivariant, i.e., the homomorphism (\ref{rho})
lifts into a homomorphism $\tilde\r:\C^\times\ra\Aut(p_*\cL\ra\ti\cX\ra\C)$,
and all the diagrams commute.

\v
Let
$D=\{w\in\C: |w|\leq 1\}$, $D^\times=\{w\in D: w\not=0\}$,
\be
M=X\times D^\times,
\ee
$\pi_{D^\times}: M\ra D^\times$ the
projection onto the second component,
$\pi_{X}: M\ra X$ the
projection onto the second component,
 and $T$ a test configuration as above.
Let $\cX_D=\pi^{-1}(D)$ and $\cX^\times_D=\pi^{-1}(D^\times)$.
We observe
there is a biholomorphic map
$\kappa:  M  \ra \cX_D^\times$ given by $(z,w)\ra \r(w)z$. The
map $(l,w)\mapsto \r(w)l$ defines
defines an isomorphism $\pi_X^*L=\kappa^*(\cL|_{\cX^\times})$. Thus
we may view $M$ as an open subset of $\cX_D$ . Moreover, if
$p:\ti\cX_D\ra\cX_D\ra D$ is a resolution of singularities,
then $\ti\cX_D^\times=p^{-1}\cX_D^\times\ra \cX_D^\times$ is biholomorphic, so
we may view $M$ as an open subset of $\ti\cX_D$ as well.
In order to avoid cumbersome notation, we shall write
$\cX,\ti\cX$ for $\cX_D$ and $\ti\cX_D$ when there is no fear of confusion.

\v

Now fix $h_0\in\ti\cH$ and let $\o_0\in \cH$ be its curvature. Let
$M=X\times D^\times$, $\pi_M: M\ra D^\times$ the
projection onto the second component and
\be
\O_0=\pi_M^*\o_0.
\ee

\begin{theorem}
\label{one} There exists a unique $C^{1,1}$ function
$\Psi: M\ra \R$ satisfying:

\v

{\rm(a)} The current $\O_0+{\sqrt{-1}\over 2}\ddb\Psi$ is non-negative on $M$.

\v
{\rm(b)} It solves the following Dirichlet problem for the completely
degenerate Monge-Amp\`ere equation on $M$
\be\label{eq}
(\O_0+{\sqrt{-1}\over 2}\ddb\Psi)^{n+1}=0\ ;\ \ \ \Psi|_{X\times\pl D}=0
\ee

\v
{\rm(c)} The current $\O_0+{\sqrt{-1}\over 2}\ddb\Psi$
is the restriction to $\ti\cX^\times$
of a  non-negative current $\Omega_{\ti\cX}$,
which is a solution of the equation
\be
\O_{\ti\cX}^{n+1}=0
\ee
on the smooth
manifold~$\ti\cX$. Moreover,
$\O_{\ti\cX}=\O'+\ddb\Psi'$
where $\O'$ is a smooth K\"ahler
metric on $\ti\cX$ and $\Psi':\ti\cX\ra\R$
is a $C^{1,1}$ function. 

\v

{\rm(d)} The function $\Psi$ is invariant under the
rotations on $D$. In particular, let $\phi_t=\Psi(z,e^{-t})$, $t\in[0,\infty)$.
Then
$\ti h_t=h_0e^{-\phi_t}$ is a $C^{1,1}$
geodesic ray in $\ti\cH$ which emanates from $h_0$, i.e.,
$\phi_t$ satisfies the geodesic equation
\be
\label{geodesic}
\ddot\phi_t-g_{\phi_t}^{j\bar k}\dot\phi_j\dot\phi_{\bar k}=0
\ \ {\rm on}\ X\times[0,\infty).
\ee
where $g_{\phi_t}$ is the K\"ahler metric $-\ddb\log \ti h_t$.
If $T$ is a non-trivial test configuration, then
the geodesic ray is infinite.
\end{theorem}
\v
%Let $\ti\pi:\ti\cX\ra\C$ be the map $\ti\pi=\pi\circ p$ and let
% $\ti\cX^\times=\ti\pi^{-1}(D^\times)$.

\v

\section{Proof of theorem}
\setcounter{equation}{0}

We divide the proof of the theorem into several steps. The first step
is the following lemma, which is an extension of the results
of Chen \cite{Ch} to the case of general complex manifolds with boundary:

\begin{lemma}
\label{monge}
Let $\ti\cX$ be a compact complex manifold with smooth boundary $\partial\ti\cX$.
If $\ti\cX$ admits a smooth K\"ahler metric $\Omega$,
then the Dirichlet problem
\be
\label{dirichlet}
(\Omega+{\sqrt{-1}\over 2}\ddb\Phi)^{n+1}=0 \ {\rm on}\ \ti\cX,\ \ \
\Phi_{\big\vert_{\partial\ti\cX}}=0
\ee
admits a unique $C^{1,1}$ solution $\Phi$ which
is $\Omega$-plurisubharmonic, i.e., $\Omega+{\sqrt{-1}\over 2}\ddb\Phi\geq 0$.
\end{lemma}

{\it Proof of Lemma \ref{monge}}:
Following Chen \cite{Ch}, we let
$t\in [0,1]$ and consider the equation

\be\label{t}  ({\Omega}+{\sqrt{-1}\over 2}\ddb\Phi_t)^{n+1}\ = \ t\cdot {\O}^n\ ; \ \
\Phi_t=0\ \ {\rm on}
\  \ \pl\ti\cX
\ee

Then (\ref{t})
has a solution when $t$=1, and, since $\Omega$ is a positive definite $(1,1)$-form,
the equation is elliptic and the set of $t$ for which
a solution exists is open. Thus, to prove the theorem, it suffices to show that
that $||\Phi_t||_{C^2}$ is uniformly bounded in $t$.
This follows in turn from the estimates of \cite{CKNS}, \cite{G}, \cite{Y78} and \cite{Ch}.
The uniqueness statement is a general property of the Dirichlet
problem for complex Monge-Amp\`ere equations (see e.g. \cite{PS06}, Theorem 6).
Q.E.D.

\bigskip
We note that the condition that $\Omega$ is a K\"ahler metric implies
the existence of a subsolution $\underline{\Phi}$ of the Dirichlet
problem,
\be
\label{subsolution}
(\Omega+{\sqrt{-1}\over 2}\ddb\underline{\Phi})^{n+1}>0 \ {\rm on}\ \cM,\ \ \
\underline{\Phi}_{\big\vert_{\partial\cM}}=0
\ee
namely $\underline{\Phi}=0$. From this point of view,
the formulation of
Lemma \ref{monge} is modeled on the earlier result of Guan \cite{G}, reducing
the solution of the complex Monge-Amp\`ere equation
to finding a subsolution.

\bigskip
We shall apply Lemma \ref{monge} to the case when $\ti\cX$ is the
equivariant resolution of the test configuration $T$
introduced in the previous section, with boundary
$\partial\ti\cX=(p\pi)^{-1}(\partial D)$. Thus
we need to construct a K\"ahler metric $\Omega$ on $\ti\cX$,
in fact a metric with some precise boundary conditions.
To do this, we prove next
the following lemma:

\medskip

\begin{lemma}
\label{bu}
There exists a line bundle $\cM\ra\ti\cX$ and an integer $m>0$ with the following properties:
\v
(i) $p^*\cL^m\otimes \cM\ra\ti\cX$ is ample, that is, it has a metric
of positive curvature.
\v
(ii) $\cM|_{\ti\cX^\times}$ is trivial, that is, there is a meromorphic
section $\m:\ti\cX\ra \cM$ whose restiction to $\ti\cX^\times$
is holomorphic and nowhere vanishing.
\end{lemma}

\bigskip

Before establishing the lemma, we recall the definition of a blow-up.
If $W$ is a scheme, and if $Z\sub W$ is a closed subscheme with ideal
sheaf $I_Z\sub \cO_W$, then $\cB(W;Z)$,
the blow-up of $W$ along $Z$, is the
scheme
${\rm Proj}\left(\oplus_{d=0}^\i
I_Z^d\right)$. If $W$ and $Z$ are
both smooth varieties, then $\cB(W;Z)$
is also a smooth variety whose underlying set is the
disjoint union
$\cB(W;Z)=(W\backslash Z)\cup \P(E)$.
Here $E\ra Z$ is the vector bundle
$E=TW|_Z/TZ$. If we let $p:\cB(W;Z)\ra W$
be the projection map, then $p$ is
biholomorphic over the set $W\backslash Z$.
\v
In order to describe the complex
structure in a neighborhood of points
in $p^{-1}(Z)$, let $a\in Z$ and let
$U\sub W$ be a small open set containing
$a$.
Choose local coordinates $(z_1,...,z_n)$
for $W$ centered at $a$ with the property
$Z\cap U=\{z_{1}=\cdots =z_{n-d}=0\}$. Let
$B_U=\{ (z,y)\in U\times \P^{n-d-1}:
y_iz_j-y_jz_i=0$ for all $1\leq i,j\leq n-d\}$. Then $B_U$  is a submanifold of
$U\times \P^{n-d-1}$. Define
$\phi_U:B_U\ra \cB(W;Z)$ as follows: $\phi_U(z,y)=z $ if $z\notin Z$
and $\phi_U(z,y)=\s_{j=1}^{n-d} y_j{\pl\over \pl z_j}$ if $z\in Z$. Then $\phi_U$
maps $B_U$ biholomorphically onto
the open set $p^{-1}(U)\sub \cB(W;Z)$.
This defines the complex structure
on $\cB(W;Z)$.

\bigskip
{\it Proof of Lemma \ref{bu}.} The map $p:\ti\cX\ra\cX$ is a composition of blow-ups
with smooth centers. Thus, by induction, we may assume $p$ is a single
blow-up: Let $W\sub\P^N$ be a smooth variety and $Z\sub W$ a smooth
subvariety. Let $L=O(1)|_W$ and $p: W'=Bl(W;Z)\ra W$ the blow-up of
$W$ with center $Z$. We wish to show that there exists a line bundle
$M\ra W'$ such that $p^*L^m\otimes M\ra W'$ is ample and $M_{W'\backslash E}$
is trivial, where $E\sub W'$ is the exceptional divisor.
Let $p_1: Bl(\P^N,Z)\ra \P^N$ be the blow up of $\P^N$ along $Z$.
Since $Bl(W,Z)\sub Bl(\P^N,Z)$ is the closure of $p_1^{-1}(W\backslash Z)\sub Bl(\P^N,Z)$,
we may assume $W=\P^N$.
\v
Thus we let $Z\sub\P^N$ be a smooth subvariety of dimension $d$. Then $Z$
is a local complete intersection. This
means that there exists an integer $k>0$, polynomials $F_1,....,F_r$ of
degree $k$, and polynomials $f_1,....,f_m$ of degree $k$, with the following
properties: The $F_s$  have no common zeros and for every $s$ with $1\leq s\leq r$,
there exists $1\leq i_1\leq \cdots \leq i_{N-d}$ such that $f_{i_1},...,f_{i_{N-d}}$
generate the ideal sheaf of $Z\cap U_s$ where $U_s\sub\P^N$ is the affine open
set defined by $F_s\not=0$. In other words, if $G$ is a homogeneous polynomial
whose degree is a multiple of $k$, and if $G$ vanishes on $Z$, then there exists
and integer $a>0$ and homogeneous polynomials $A_1,...,A_{N-d}$ such that

$$ F_s^aG\ = \ A_1f_{i_1}+\cdots+A_{N-d}f_{i_{(N-d)}}
$$

Note that
$m\geq N-d$ .
Let

\be\label{calB}
\cB\ = \ \{(x,Y)\in \P^N\times \P^{m-1}: Y_if_j(x)-Y_jf_i(x)=0, 1\leq i,j\leq m\}
\ee
If $m=N-d$ (i.e., $Z$ is a complete intersection) then, by the definition of
blow-up,  $\cB=Bl(\P^N;Z)$. In
general, let
let $\cB_0\sub \P^N\times\P^{m-1}$ be the set $\cB_0=\{(x,y)\in\cB: x\notin Z\}$.
Thus $\cB_0\ra \P^N\backslash Z$ is biholomorphic.
Let $\bar\cB_0\sub \P^N\times \P^{m-1}$ be the closure of $\cB_0$ and
$P:\bar\cB_0\ra\P^N$ projection onto the first component. We claim
that $\bar\cB_0=Bl(\P^N;Z)$. To see this, define
$\T:\cB(\P^N;Z)\backslash p^{-1}(Z)\ra \bar\cB_0\backslash P^{-1}(Z)$
as follows: $\T=P^{-1}\circ p$. We must show that $\T$ uniquely extends
to a biregular map
$\T:\cB(\P^N;Z)\ra \bar\cB_0$
\v

To prove that $\T$ has a unique extension, we may work locally:
Fix $z\in Z$. Choose $F_s$ , such
that $F_s(z)\not=0$, and choose $i_1<\cdots  <i_{N-d}$ as above. Without
loss of generality, we may assume $i_j=j$.  We shall write $F=F_s$ and $U=U_s$.
Then

\be p^{-1}(U\backslash Z)\ = \ \{(x,y)\in (U\backslash Z)\times \P^{N-d-1}:
y_if_j(x)-y_jf_i(x)=0, \ 1\leq i,j\leq N-d\ \}
\ee
\be P^{-1}(U\backslash Z)=\{(x,Y)\in (U\backslash Z)\times\P^{m-1}:
Y_if_j(x)-Y_jf_i(x)=0, 1\leq i,j\leq m\}
\ee
We must examine the map $\T: p^{-1}(U\backslash Z)\ra  P^{-1}(U\backslash Z)$.
If $\T(x,y)=(x',Y)$ then $x'=x$  and $Y_i=y_i$ if $1\leq i\leq N-d$. Moreover,
on the open set $y_l\not=0$ ($1\leq l\leq N-d)$, we have

\be\label{ext}
 Y_i\ = \ {y_lf_i(x)\over f_l(x)}\ = \ {y_l\s_{i=1}^{N-d} A_i(x)f_i(x)\over F^s(x)f_l(x)}
\ = \ {\s_{i=1}^{N-d} A_i(x)y_i(x)\over F^s(x)}
\ee
Note that $f_l(x)\not=0$ in (\ref{ext}) since $y_if_l=y_lf_i$ so if $f_l(x)=0$
then $f_i(x)=0$ for all $1\leq i\leq N-d$ which implies $x\in Z$, and contradiction.
\v
Thus (\ref{ext}) gives the formula for $\T(x,y)$ when
$x\notin Z$. On the other hand, the denominator $F^s(x)$
does not vanish if $x\in Z$ and thus (\ref{ext}) defines
a unique extension of $\T$ from $p^{-1}(U\backslash Z)$ to $p^{-1}(U)$.
The map $\T$ is clearly biregular. In fact, its inverse is given
explicitly by $(x',Y)\ra (x,y)$
where $x=x'$ and $y_i=Y_i$ for $1\leq i\leq N-d$.

\v

Now we can finish the proof of Lemma \ref{bu}:  We must show that there
exists a line bundle $\cM\ra \cB(\P^N;Z)$ such that
$p^*O_{\P^N}(m)\otimes \cM\ra \cB(\P^N;Z)$
is ample, for some $m>0$, and such that $\cM|_{\cB(\P^N;Z)\backslash p^{-1}(Z)}$
is trivial. Since $\cB=\cB(\P^N;Z)=\bar\cB_0\sub\P^N\times\P^{m-1}$, we see
that $p^*O_{\P^N}(1)\otimes q^*O_{\P^{m-1}}(1)\ra \cB$ is ample, where
$q:\cB\ra\P^{m-1}$ the the projection onto the second factor.
Let $\cM^{-1}=p^*O_{\P^N}(k)\otimes q^*O_{\P^{m-1}}(-1)$ and let
$s$ be the section of
$\cM^{-1}$ defined by $s={f_j(x)\over y_j}$ on the open set
$y_j\not=0$, then $s$ is a global section which vanishes precisely on
the exceptional divisor $D=p^{-1}(Z)\sub\cB$. Thus $\cM$ is trivial
on the complement of the exceptional divisor. Moreover

\be p^*O_{\P^N}(k+1)\otimes M\ = \ p^*O_{\P^N}(1)\otimes q^*O_{\P^{m-1}}(1)
\ee

This proves  Lemma \ref{bu}.

\v

Now let $T$ and $M$ be as in Lemma \ref{bu}, and let $\o_0\in c_1(L)$
be a fixed K\"ahler metric and let $h_0\in\ti\cH$ be a hermitian
metric on $L$ such that $\o_0$ is the curvature of $h_0$.

\v

\begin{lemma}\label{sub}
There exists a K\"ahler metric $\O$ on $\ti\cX$ with the following
properties:
\v
(i) $m\O$ is the curvature of a hermitian metric $H$ on $p^*\cL^m\otimes \cM$.

\v

(ii) $\r(w)^*(\O|_{X_w})=\o_0$ for all $w\in \pl D$.

\v
(iii) $\O$ is invariant under the group action.

\end{lemma}

{\it Proof.} According to Lemma \ref{bu}, there is a metric $H_1$ on
$p^*\cL^m\otimes \cM$ with positive curvature $\O_1$. We shall modify $H_1$
in a neighborhood of $\pl\ti\cX$ in order to produce a new positively
curved metric whose boundary value is $\o_0$: Let $s$ be a meromorphic
section of $\cM$ which is holomorphic and nowhere vanishing on $\ti\cX^\times$.
We define a metric $H_2$ on $\pi_X^*L^m=L^m\times D^\times$ as follows:
\be H_2(\ell,w)\ = \ H_1(p^*(\r(w)\ell)\otimes s))
\ee
Since $H_2$ and $h_0^m$ are two metrics on the same line bundle
$L^m\times D^\times$, there exists a smooth function $\Psi: M=X\times D^\times\ra\R$
such that $H_2=h_0^me^{-\Psi(x,w)}$.
\v
Let $\eta: D^\times\ra [0,1]$ be a smooth function such that $\eta(w)=1$
if $|w|\leq {1\over  3}$ and $\eta(w)=0$ if $|w|\geq {2\over 3}$.
Let
\be
H_3=h_0^me^{-\eta(w)\Psi(x,w)}
\ee
Then $H_3=h_0$ if $|w|\geq {2\over 3}$ and $H_3=H_2$ if $|w|\leq {1\over 3}$.
Moreover, the curvature of $H_3$ is positive on the fibers $X\times \{w\}$
for all $w\in D^\times$.
\v
Let $\al>0$ be a large positive number, and define

$$ H_4 \ = \ H_3e^{-\al(|w|^2-1)}
$$
Then $H_4(\ell,w)\ = \ H_5(p^*(\r(w)\ell)\otimes s))$ for some
unique smooth metric $H_5$ on $p^*\cL^m\otimes \cM$, and $H$ satisfies properties
{\it (i)} and {\it (ii)} by construction.
Finally, we average $H_5$ over the $S^1$ action to obtain the desired metric $H$. Q.E.D.

\bigskip
{\it Proof of theorem.}
We now turn to the proof of Theorem 1: Fix $\O$ as in Lemma \ref{sub}.
By Lemma \ref{monge}, there is a unique $\Omega$-plurisubharmonic,
$C^{1,1}$ invariant
solution $\Phi$ on $\ti\cX$
of the equation (\ref{dirichlet}). Now, restricted to $\ti\cX^\times\sim X\times D^\times$,
the metric $\O$ can be expressed as
\be
\O=\O_0+{\sqrt{-1}\over 2}\ddb \Phi_0
\ee
for some smooth function $\Phi_0\in C^\infty(X\times \bar D^\times)$, with
$\Phi_0=0$ on $X\times \partial D$.
The equation (\ref{dirichlet}) on $\ti\cX^\times$ can be written as
\be
(\O_0+{\sqrt{-1}\over 2}\ddb (\Phi_0+\Phi))^{n+1}=0,
\ \ \
\ee
Set $\Psi=\Phi_0+\Phi$. Then $\Psi$ satisfies conditions (a) and (b)
of the theorem. The uniqueness stated in (c) is a consequence of
the uniqueness of solutions of the Monge-Amp\`ere equation
on $\ti\cX$ stated in Lemma \ref{monge}. Clearly, the invariance of
$\Omega$ and of $\Phi$ results in the invariance of $\Psi$.
It is now a well-known fact that, for invariant functions $\Psi$
on $X\times D$, the Monge-Amp\`ere equation in (\ref{eq})
reduces precisely to the geodesic equation (\ref{geodesic}). Q.E.D.

$^*$ Department of Mathematics

Columbia University, New York, NY 10027

\v

$^{\dagger}$ Department of Mathematics

Rutgers University, Newark, NJ 07102

\enddocument


\newpage

\begin{thebibliography}{99}

{\small

\bibitem{AT} Arezzo, C., and Tian, G.,
``Infinite geodesic rays in the space of K\"ahler potentials'',
Ann. Sci. Norm. Sup. Pisa Sci. (5) 2 (2003) 617-630,
arXiv : math.DG / 0210389.

\bibitem{BM} Bierstone, E. and P. Milman, ``Canonical desingularization in
characteristic zero by blowing up the
maximal stata of a local invariant. Inv. Math. {\bf 128} (1997), 207-302
\bibitem{BT76} Bedford, E. and B.A. Taylor,
``The Dirichlet problem for a complex Monge-Amp\`ere equation'',
Invent. Math. {\bf 37} (1976), 1-44.

\bibitem{BT82} Bedford, E. and B.A. Taylor,
``A new capacity for plurisubharmonic functions'', Acta Math. {\bf 149}
(1982), 1-40.



\bibitem{Ch} Chen, X.X.,
``The space of K\"ahler metrics'', J. Differential Geom.
{\bf 56} (2000), 189-234.

\bibitem{Ch06} Chen, X.X.,
``Space of K\"ahler metrics III: On the lower bound
of the Calabi energy and geodesic distance'',
arXiv: math.DG / 0606228

\bibitem{CKNS} Caffarelli, L., J.J. Kohn, L. Nirenberg, and J. Spruck,
``The Dirichlet problem for non-linear second-order elliptic equations II.
Complex Monge-Amp\`ere and uniformly elliptic equations'',
Comm. Pure Appl. Math. {\bf 38} (1985) 209-252.


\bibitem{D99} Donaldson, S.K.,
``Symmetric spaces, K\"ahler geometry, and Hamiltonian
dynamics'', Amer. Math. Soc. Transl. {\bf 196} (1999) 13-33.


\bibitem{D02} Donaldson, S.K.,
``Scalar curvature and stability of toric varieties'',
J. Differential Geom. {\bf 62} (2002) 289-349


\bibitem{G} Guan, B.,
 ``The Dirichlet problem for the complex Monge-Amp\`ere operator'',
Comm. Anal. Geom. {\bf 6} (1998) 687-703

\bibitem{H} Hauser, H., ``The Hironaka theorem on resolution of singularities'',
Bull. Am. Math. Soc. {\bf 40} (2003)
 323-403

\bibitem{M87} Mabuchi, T.,
``Some symplectic geometry on compact K\"ahler manifolds I'',
Osaka J. Math. {\bf 24} (1987) 227-252.


\bibitem{PS06} Phong, D.H. and J. Sturm,
``The Monge-Amp\`ere operator and geodesics in the space of K\"ahler potentials'',
Inventiones Math. {\bf 166} (2006) 125-149,
arXiv: math.DG / 0504157.

\bibitem{PS06a} Phong, D.H. and J. Sturm,
``Test Configurations for K-Stability and Geodesic Rays''
preprint, math.DG/0606423


\bibitem{S92} Semmes, S.,
``Complex Monge-Amp\`ere equations and symplectic manifolds'',
Amer. J. Math. {\bf 114} (1992) 495-550.

\bibitem{SZ} Song, J. and S. Zelditch,
``Bergman metrics and geodesics in the space of K\"ahler metrics
on toric varieties'', arXiv:0707.3082 (math.CV)


\bibitem{SZf} Song, J. and S. Zelditch,
``Test configurations and geodesic rays on toric varieties'',
in preparation.



\bibitem{T90} Tian, G.,
``On a set of polarized K\"ahler metrics on algebraic manifolds'',
J. Differential Geom. {\bf 32} (1990) 99-130


\bibitem{T97} Tian, G.,
``K\"ahler-Einstein metrics with positive scalar curvature'',
Inventiones Math. {\bf 130} (1997) 1-39



\bibitem{V} Villamayor, O., ``Constructiveness of Hironaka's resolution'',
Ann. Scient. Ec. Norm. Sup. Paris {\bf 22}
(1989) 1-32

\bibitem{Y78} Yau, S.T.,
``On the Ricci curvature of a compact K\"ahler manifold
and the complex Monge-Amp\`ere equation I'',
Comm. Pure Appl. Math. {\bf 31} (1978) 339-411.

\bibitem{Y} Yau, S.T., ``Open problems in geometry'',
Proc. Symp. Pure Math. {\bf 54}, Amer. Math. Soc.
Providence, RI (1993) 1-28


\bibitem{Z} Zelditch, S.,
``The Szeg\"o kernel and a theorem of Tian'',
Int. Math. Res. Notices {\bf 6} (1998) 317-331.

}


\end{thebibliography}
\enddocument